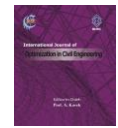

# TOWER CRANES AND SUPPLY POINTS LOCATING PROBLEM USING CBO, ECBO, AND VPS

A. Kaveh[*, †] and Y. Vazirinia
*Centre of Excellence for Fundamental Studies in Structural Engineering, Iran University of Science and Technology, Narmak, Tehran, P.O. Box 16846-13114, Iran*

## ABSTRACT

Tower cranes are major and expensive equipment that are extensively used at building construction projects and harbors for lifting heavy objects to demand points. The tower crane locating problem to position a tower crane and supply points in a building construction site for supplying all requests in minimum time, has been raised from more than twenty years ago. This problem has already been solved by linear programming, but meta-heuristic methods spend less time to solving the problem. Hence, in this paper three newly developed meta-heuristic algorithms called CBO, ECBO, and VPS have been used to solve the tower crane locating problem. Three scenarios are studied to show the applicability and performance of these meta-heuristics.

## 1. INTRODUCTION

Every construction project needs enough spaces for temporary facilities in order to perform the construction activities in a safe and efficient manner. Construction site-level facilities layout is an important step in site planning. Planning construction site spaces to allow for safe and efficient working status is a complex and multi-disciplinary task as it involves accounting for a wide range of scenarios. Construction site layout problems are known as combinatorial optimization problems. There are two types of procedures for solution consisting of meta-heuristics for large search sized problems and the exact method with global search for smaller search sized problems [1]. For example, Li and Love [2] developed a construction site-level facility layout problem for allocating a a set of predetermined facilities into a set of predetermined locations, while satisfying the layout constraints and requirements. They used a genetic algorithm to solve the problem by assuming that the predetermined locations are in rectangular shape and are large enough to accommodate the largest facility. Gharaie et al. [3] resolved their model by Ant Colony Optimization, and Kaveh et al. [4] used Colliding Bodies Optimization and its enhanced version. Similarly, Cheung et al. [5] have developed another model for construction site layout planning and solved it by Genetic Algorithm. Also Liang and Chao [6], Wong et al. [7], Kaveh et al. [8], and Kaveh et al. [4] have employed Multi-searching tabu search, Mixed Integer Programming, and CBO, ECBO, and PSO, respectively.

---
[*]Corresponding author: Centre of Excellence for Fundamental Studies in Structural Engineering, Iran University of Science and Technology, Narmak, Tehran, P.O. Box 16846-13114, Iran
[†]E-mail address: alikaveh@iust.ac.ir (A. Kaveh)

The tower crane is an important facility used in the vertical and horizontal transportation of materials, particularly heavy prefabrication units such as steel beams, ready-mixed concrete, prefabricated elements, and large-panel formworks [9]. The Tower Crane Locating Problem to positioning a tower crane and supply points in a building construction site for supplying all the requests in a minimum time, has been raised more than twenty years ago. An analytical model was developed by Zhang et al. [10] considering the travel time of tower crane hooks and adopting a Monte Carlo simulation to optimize the tower crane location. However, considered tower crane in their study was a single one and the effect of location of supply points on lifting requirements and travel time has been neglected. Tam et al. [9] employed an artificial neural network model for predicting tower crane operations and next they used a genetic algorithm model to optimize the crane and supply points layout [9] and [11]. The case study used by Tam et al. [11] to show the effectiveness of their model was subsequently used in a number of researches to compare the effectiveness of other optimization methods. For example, Huang et al. [1] used mixed integer linear programming (MILP) to optimize the crane and supply locations and showed that their method reduced the travel time of the hook by 7% compared to the results obtained from the previous genetic algorithm. MILP was used to ensure achieving a global optimal solution [1]. In this research, we use new meta-heuristic techniques, because they require less computational time to solving the problem.

Solving real-life problems by meta-heuristic algorithms has become to an interesting topic in recent years. Many meta-heuristics with different philosophy and characteristics are developed and applied to a wide range of fields. The aim of these optimization methods is to efficiently explore the search space in order to find global or near-global solutions. Since these algorithms are not problem specific and do not require the derivatives of the objective function, they have received increasing attention from both academia and industry [12]. Meta-heuristic methods are global optimization methods that try to reproduce natural phenomena (Genetic Algorithm [13], Particle Swarm Optimization [14], Water Evaporation Optimization [15]), humans social behavior (Imperialist Competitive Algorithm[16]), or physical phenomena (Charged System Search (CSS) [17], Colliding Bodies Optimization [18], Big Bang-Big Chrunch [19], Vibrating Particle System (VPS) [12]). Exploitation and exploration are two important characteristics of meta-heuristic optimization methods, Kaveh [20]. Exploitation serves to search around the current best solutions and to select the best possible points, and Exploration allows the optimizer to explore the search space more efficiently, often by randomization.

In this paper three newly developed meta-heuristics named: Colliding Bodies Optimization (CBO), Enhanced Colliding Bodies Optimization (ECBO) [21], and Vibrating

Particle System (VPS) [12] are used for tower crane and material supply locating model that proposed by Huang et al. [1] and results are discussed.

## 2. OPTIMIZATION ALGORITHMS

### 2.1 Colliding bodies optimization

An efficient algorithm, inspired from the momentum, and energy rules of the physics, named Colliding Bodies Optimization, that was developed by Kaveh and Mahdavi [18]. CBO does not depend on any internal parameter and also it is extremely simple in the sense. In this method, one body collides by another body and they moves to the lower cost. Each solution candidate "X" at CBO, contains a number of variables (i.e., $X_i = \{X_{i,j}\}$) and is considered as a colliding body (CB). The massed bodies are divided in two main equal groups; i.e., stationary and moving bodies (Fig. 2), where the moving bodies moves to stationary bodies and a collision occurs between the pairs of bodies. The goal of this process is: (i) to improve the locations of moving bodies and (ii) to push stationary bodies toward the better locations. After the collision, new locations of colliding bodies are updated based on the new velocity by using the collision rules. The main procedure of the CBO is described as:

Step 1: The initial positions of colliding bodies are determined with random initialization of a population of individuals in the search space:

$$x_i^o = x_{min} + rand \times (x_{max} - x_{min}), i = 1,2,3, \dots, n \qquad (1)$$

where, $x_i^o$ determines the initial value vector of the $i$th colliding body. $x_{max}$ and $x_{min}$ are the minimum and the maximum allowable values vectors of variables, respectively; rand is a random number in the interval [0,1]; and n is the number of colliding bodies.

Step 2: The magnitude of the body mass for each colliding body is defined as:

$$m_k = (1/_{fit(i)})/(1/\sum_{i=1}^{n}(1/_{fit(i)})) \qquad (2)$$

where fit (i) represents the objective function value of the colliding body i; n is the population size. It seems that a colliding body with good values exerts a larger mass than the bad ones. Also, for maximization, the objective function fit(i) will be replaced by $1/_{fit(i)}$.

Step 3: Then colliding bodies objective function values are arranged in an ascending order. The sorted colliding bodies are divided into two equal groups:

• The lower half of the CBs (stationary CBs); These CBs are good agents which are stationary and the velocity of these bodies before collision is zero. Thus:

$$v_i = 0, i = 1,2,3, \dots, \frac{n}{2} \qquad (3)$$

• The upper half of CBs (moving CBs): These CBs move toward the lower half. Then, according to Fig. 1, the better and worse CBs, i.e. agents with upper fitness value, of each group will collide together. The change of the body position represents the velocity of these bodies before collision as:

$$v_i = x_i - x_{i-\frac{n}{2}} , i = \frac{n}{2} + 1, \dots , n \qquad (4)$$

where, $v_i$ and $x_i$ are the velocity and position vector of the $i$th CB in this group, respectively; $x_{i-\frac{n}{2}}$ is the $i$th CB pair position of $x_i$ in the previous group.

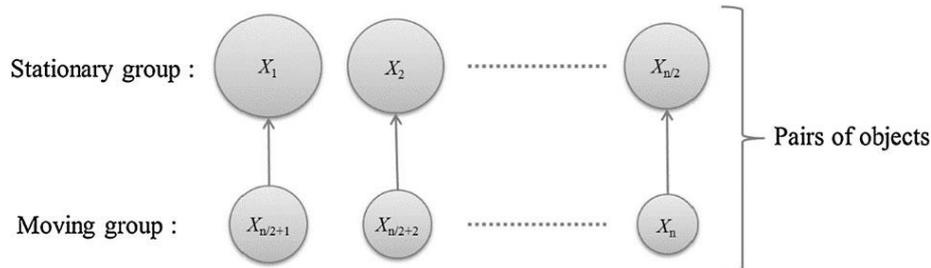

Figure 1. Pairs of CBs for collision

Step 4: After the collision, the velocities of the CBs in each group are evaluated as:
• Stationary CBs:

$$v_i' = \frac{\left(m_{i+\frac{n}{2}} - \varepsilon m_{i+\frac{n}{2}}\right) \times v_{i+\frac{n}{2}}}{\left(m_i + m_{i+\frac{n}{2}}\right)} , i = 1,2,3, \dots , \frac{n}{2} ; \qquad (5)$$

• Moving CBs:

$$v_i' = \frac{\left(m_i - \varepsilon m_{i-\frac{n}{2}}\right) \times v_i}{\left(m_i + m_{i-\frac{n}{2}}\right)} , i = \frac{n}{2} + 1, \frac{n}{2} + 2, \frac{n}{2} + 3, \dots , n ; \qquad (6)$$

where $\varepsilon$ is the Coefficient Of Restitution (COR) of the two colliding bodies, defined as:

$$\varepsilon = 1 - \frac{iter}{iter_{max}} \qquad (7)$$

with iter and $iter_{max}$ being the current iteration number and the total number of iteration for optimization process, respectively.

New positions of CBs are updated using the generated velocities after the collision in position of stationary CBs, as follow for each group:
• Moving CB:

$$x_i^{new} = x_{i-\frac{n}{2}} + rand \circ v_i' \,, i = \frac{n}{2} + 1, \frac{n}{2} + 2, \frac{n}{2} + 3, \dots, n; \qquad (8)$$

where, $x_i^{new}$ and $v_i'$ are the new position and the velocity after the collision of the $i$th moving CB, respectively; $x_{i-\frac{n}{2}}$ is the old position of the $i$th stationary CB pair.

- Stationary CB:

$$x_i^{new} = x_i + rand \circ v_i' \,, i = 1,2,3, \dots, \frac{n}{2}; \qquad (9)$$

where, $x_i^{new}$, $x_i$ and $v_i'$ are the new positions, previous positions and the velocity after the collision of the $i$th CB, respectively. rand is a random vector uniformly distributed in the range of [-1,1] and the sign '$\circ$' denotes an element-by-element multiplication.

Step 6: The process is repeated from step 2 until one termination criterion is satisfied. Termination criterion is the predefined maximum number of iterations. After getting the near-global optimal solution, it is recorded to generate the output.

---

**Pseudo Code of Colliding Bodies Optimization**

    Initial location is created randomly by Eq. (1)
    The value of the objective function is evaluated and masses are defined by Eq. (2)
    **While** stop criteria is not attained (like max iteration)
        **for** each CBs
            Calculate Stationary and moving CBs velocity before collision according Eqs. (3) and (4)
            Calculate CBs velocity after collision according by Eqs. (5) and (6)
            Update CBs position according Eqs. (8) and (9)
        **End** for
    **End** while
se 1.**End**

---

Figure 2. Pseudo code of the colliding bodies optimization

## 2.2. Enhanced colliding bodies optimization

In order to improve CBO to get faster and more reliable solutions, Enhanced Colliding Bodies Optimization (ECBO) was developed which uses memory to save a number of historically best CBs and also utilizes a mechanism to escape from local optima [11]. The pseudo of ECBO is shown in Fig. 3 and the steps involved are given as follows:

Step 1: Initialization

Initial positions of all CBs are determined randomly in an m-dimensional search space by Eq. (1).

Step 2: Defining mass

The value of mass for each CB is evaluated according to Eq. (2).

Step 3: Saving

Considering a memory which saves some historically best CB vectors and their related mass and objective function values can improve the algorithm performance without increasing the computational cost [16]. For that purpose, a Colliding Memory (CM) is utilized to save a number of the best-so-far solutions. Therefore in this step, the solution vectors saved in CM are added to the population, and the same numbers of current worst CBs are deleted. Finally, CBs are sorted according to their masses in a decreasing order.

Step 4: Creating groups

CBs are divided into two equal groups: (i) stationary group and (ii) moving group. The pairs of CBs are defined according to Fig. 1.

Step 5: Criteria before the collision

The velocity of stationary bodies before collision is zero (Eq. (3)). Moving objects move toward stationary objects and their velocities before collision are calculated by Eq. (4).

Step 6: Criteria after the collision

The velocities of stationary and moving bodies are calculated using Eqs. (5) and (6), respectively.

Step 7: Updating CBs

The new position of each CB is calculated by Eqs. (8) and (9).

---

**Pseudo Code of Enhanced Colliding Bodies Optimization**

---

    Initial location is created randomly by Eq. (1)

    The value of objective function is evaluated and masses are defined by Eq. (2)

    **While** stop criteria is not attained (like max iteration)

        **for** each CBs

            Calculate Stationary and moving CBs velocity before collision according Eqs. (3) and (4)

            Calculate CBs velocity after collision according by Eqs. (5) and (6)

            Update CBs position according Eqs. (8) and (9)

            **If** rand i < Pro

                One dimension of the $i$th CB is selected randomly and regenerate by Eq. (10)

            **End** if

        **End** for

    **End** while

**End**

---

Figure 3. Pseudo code of the enhanced colliding bodies optimization [24]

Step 8: Escape from local optima

Meta-heuristic algorithms should have the ability to escape from the trap when agents get close to a local optimum. In ECBO, a parameter like Pro within (0, 1) is introduced and it is specified whether a component of each CB must be changed or not. For each colliding body Pro is compared with $rn_i$ ($i = 1, 2, \ldots, n$) which is a random number uniformly distributed within (0, 1). If $rn_i < Pro$, one dimension of the ith CB is selected randomly and its value is regenerated as follows:

$$x_{ij} = x_{j,min} + rand \times (x_{j,max} - x_{j,min}), i = 1,2,3, \ldots, n; \qquad (10)$$

where $x_{ij}$ is the $j$th variable of the $i$th CB. $x_{j,min}$ and $x_{j,max}$ respectively, are the lower and upper bounds of the $j$th variable. In order to protect the structures of CBs, only one dimension is changed. This mechanism provides opportunities for the CBs to move all over the search space thus providing better diversity.

Step 9: Terminating condition check

The optimization process is terminated after a fixed number of iterations. If this criterion is not satisfied go to Step 2 for a new round of iteration.

### 2.3 Vibrating particle system

The VPS is a population-based algorithm which simulates a free vibration of single degree of freedom systems with viscous damping [12]. Similar to other multi-agent methods, VPS has a number of individuals (or particles) consisting of the variables of the problem. In the VPS each solution candidate is defined as "X", and contains a number of variables (i.e., $X_i = \{X_i^j\}$) and is considered as a particle. Particles are damped based on three equilibrium positions with different weights, and during each iteration the particle position is updated by learning from them: (i) the historically best position of the entire population (HB), (ii) a good particle (GP), and (iii) a bad particle (BP). The solution candidates gradually approach to their equilibrium positions that are achieved from current population and historically best position in order to have a proper balance between diversification and intensification. Main procedure of this algorithm is defined as:

Step 1: Initialization

Initial locations of particles are created randomly in an n-dimensional search space, by Eq. (11):

$$x_i^j = x_{min} + rand \times (x_{max} - x_{min}), \qquad i = 1,2,3, \dots, n;$$ (11)

where, $x_i^j$ is the $j$th variable of the particle $i$. $x_{max}$ and $x_{min}$ are respectively the minimum and the maximum allowable values vectors of variables. rand is a random number in the interval [0,1]; and n is the number of particles.

Step 2: Evaluation of candidate solutions

The objective function value is calculated for each particle.

Step 3: Updating the particle positions

In order to select the GP and BP for each candidate solution, the current population is sorted according to their objective function values in an increasing order, and then GP and BP are chosen randomly from the first and second half, respectively.

According to the above concepts, the particles position are updated by follow equation:

$$x_i^j = \omega_1 . [D.A. rand1 + HB^j] + \omega_2 . [D.A. rand2 + GP^j] + \omega_3 . [D.A. rand3 + BP^j]$$ (12)

where $x_i^j$ is the $j$th variable of the particle $i$. $\omega_1$, $\omega_2$, $\omega_3$, are three parameters to measure the relative importance of HB, GP and BP, respectively ($\omega_1 + \omega_2 + \omega_3 = 1$). rand1, rand2, and rand3 are random numbers uniformly distributed in the range of [0, 1] respectively. The parameter A is defined as:

$$A = [\omega_1 . (HB^j - x_i^j)] + [\omega_1 . (GP^j - x_i^j)] + [\omega_1 . (BP^j - x_i^j)]$$ (13)

Parameter D is a descending function based on the number of iterations:

$$D = (\frac{\text{iter}}{\text{iter}_{max}})^{-\alpha} \tag{14}$$

In order to have a fast convergence in the VPS, the effect of BP is sometimes considered in updating the position formula. Therefore, for each particle, a parameter like p within (0,1) is defined, and it is compared with rand (a random number uniformly distributed in the range of [0,1]) and if p < rand, then $\omega = 0$ and $\omega_2 = 1 - \omega_1$.

Three essential concepts consisting of self-adaptation, cooperation, and competition are considered in this algorithm. Particles moves towards HB so the self-adaptation is provided. Any particle has the chance to have influence on the new position of the other one, so the cooperation between the particles is supplied. Because of the p parameter, the influence of GP (good particle) is more than that of BP (bad particle), and therefore the competition is provided.

Step 4: Handling the side constraints

There is a possibility of boundary violation when a particle moves to its new position. In the proposed algorithm, for handling boundary constraints a harmony search-based approach is used [17]. In this technique, there is a possibility like harmony memory considering rate (HMCR) that specifies whether the violating component must be changed with the corresponding component of the historically best position of a random particle or it should be determined randomly in the search space. Moreover, if the component of a historically best position is selected, there is a possibility like pitch adjusting rate (PAR) that specifies whether this value should be changed with the neighboring value or not.

| Pseudo code of Vibrating Particles System (VPS) |
|---|
| Initialize algorithm parameters |
| Create initial positions randomly by Eq. (11) |
| Evaluate the values of objective function and store HB |
| **While** maximum iterations is not fulfilled |
|    **for** each particle |
|       The GP and BP are chosen |
|       **if** P<rand |
|          $W_3=0$ and $w_2=1\text{-}w1$ |
|       **end** if |
|       **for** each component |
|          New location is obtained by Eq. (12) |
|       **end** for |
|       Violated components are regenerated by harmony search-based handling approach |
|    **end** for |
|  **end** while |
| The values of objective function are evaluated and HB is updated |
| **end** |

Figure 4. Pseudo code of the vibrating particles system algorithm [22]

Step 5: Terminating condition check

Steps 2 through 4 are repeated until a termination criterion is fulfilled. Any terminating condition can be considered, and in this study the optimization process is terminated after a fixed number of iterations. Pseudo code of the VPS is showed in Fig. 4.

## 3. PROBLEM: OPTIMIZATION OF TOWER CRANE LOCATION AND MATERIAL SUPPLY POINTS

Many researches have been conducted on the locating and transporting time of a tower crane, such as: Choi and Harris [22] improved a mathematical model for determining the most suitable tower crane location; Zhang et al. [10] developed the Monte Carlo simulation approach to optimize tower crane location; Tam and Tong [9] employed an artificial neural network model for predicting tower crane operations and genetic algorithm model for site facility layout [9, 11]. Huang et al. [1] developed a mixed integer linear programming (MILP) to optimize the crane and supply locations were and their model decreased the travel time of the hook by 7% compared to the results obtained from the previous genetic algorithm. Travel distance between the supply and demand points can be calculated by the Eq. (15) through Eq. (19) referring to Figs. 5 and 6.

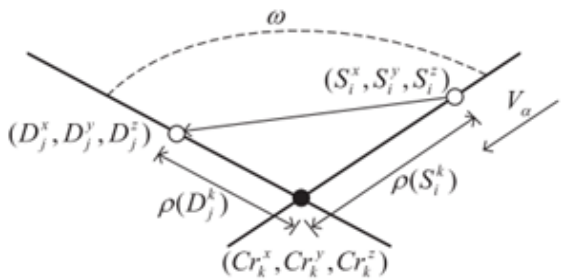

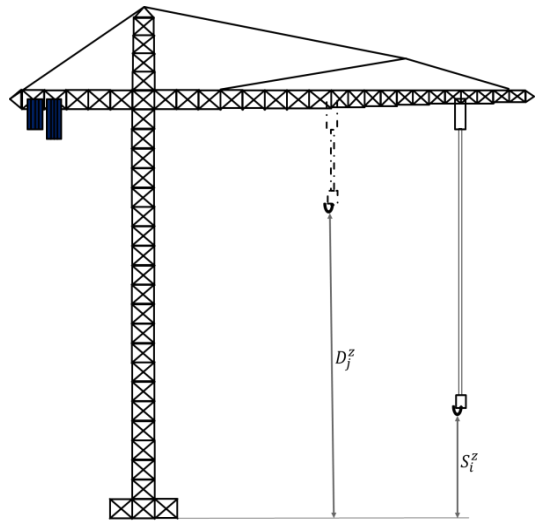

—— Rib of the crane

---- Angular movement path of the rib of the crane (tangent movement)

○  Hook position

←  Change of hook position (radial movement)

●  Crane position

Figure 5. Radial and tangent movements of the hook

Figure 6. Vertical movement of the hook

$$\rho_{(D,Cr)} = \sqrt{(X_D - X_C)^2 + (Y_D - Y_C)^2} \tag{15}$$

$$\rho_{(S,Cr)} = \sqrt{(X_S - X_C)^2 + (Y_S - Y_C)^2} \qquad (16)$$

$$L = \sqrt{(X_S - X_D)^2 + (Y_S - Y_D)^2} \qquad (17)$$

$$T_a = \frac{|\rho_D - \rho_S|}{V_a} \qquad (18)$$

$$T_\omega = \frac{1}{V_\omega} \arccos\left[\frac{L^2 + {\rho_D}^2 + {\rho_S}^2}{2 * \rho_D * \rho_S}\right], [0 < \arccos(\theta) < \pi] \qquad (19)$$

Hook movement time is an important parameter to evaluate the total time of material transportation using a tower crane. The hook movement time has split up into horizontal and vertical paths to reflect the operating costs by giving an appropriate cost-time factor. Corresponding movement paths along different directions can be seen from Figs. 5 and 6.

A continuous type parameter α indicates the degree of coordination of the hook movement in radial and tangential directions which depends on the control skills of a tower crane operator, times for horizontal and vertical hook movements can be calculated in Eqs. (20) and (21), respectively.

$$T_h = \max\{T_a + T_\omega\} + \alpha * \min\{T_a + T_\omega\} \qquad (20)$$

$$T_v = \frac{|D_D - D_S|}{V_h} \qquad (21)$$

The total travel time of tower crane at location k between supply point i and demand point j, $T_{i,j}^k$, can be calculated using Eq. (22) by specifying the continuous type parameter ß for the degree of coordination of hook movement in horizontal and vertical planes. Depending on different site conditions, skills of operators, or even the visibility level due to environmental and weathering effects, the movement of the tower crane and the hook operation may be influenced and overall efficiency can be reduced, meaning that longer operating time is required for moving the tower crane from one point to another [24]. The aggregate travel time from the material supply location to the demand point should be increased accordingly if the operator's line of sight is obstructed. To realize these site operating difficulties, another numerical parameter $\gamma_k$ is introduced to factor up the original tower crane and hook travel times given in Eq. (22). Different $\gamma_k$ may be used for different tower crane locations k to determine the location specific effects within a construction site. If advance vision system is installed in tower cranes to assist operators, the operation time can be faster and a smaller $\gamma_k$ can be set [1].

$$T_{i,j}^k = \lambda_k * [\max\{T_h + T_v\} + \beta.\min\{T_h + T_v\}] \qquad (22)$$

Huang et al. [1] provided three scenarios to demonstrate the flexibility of their proposed MILP model of a tower crane location. The formulation will be extended to consider homogeneous and non-homogeneous storage supply points where different materials can be stored in different strategies by confining the solution region with extra sets of linear type governing constraints.

Only one tower crane can be modeled which can be allocated at any one of the available locations. Binary variables like $\zeta_k$ are defined for a location k, where $\zeta_k = 1$ if the location k is selected for a tower crane location or $\zeta_k = 0$ otherwise. Constraint (23) is required, so only the best tower crane location can be picked in the optimization framework.

$$\sum_{k=1}^{K} \zeta_k = 1, \qquad \forall k \in \{1, K\} \tag{23}$$

A set of binary variables $\Delta_j$ is introduced to represent the existence of a demand location where j is the potential demand point. Depending on the input material demand profile $Q_{l,j}$ for material type l, constraint set (24) is required to ensure the binary variable $\Delta_j$ to be "1" if there is a demand at location j and "0" if the demand does not exist. 'M' is an arbitrary large integral number.

$$M\Delta_j \geq \sum_{l=1}^{L} Q_{l,j} \geq \Delta_j, \forall j \in \{1, J\} \tag{24}$$

### 3.1 Homogeneous material supply point

As a management problem, it is worth examining the total cost for transporting all the required materials to demand points through a tower crane if the materials can be stored and supplied in more than one location without setting a storage limit on various supply locations realizing that the supply locations have infinite material storage capacity, which is always the case in large scale construction sites. Under a homogeneous material supply system, each supply point provides a temporary material storage that is restricted to supply only one type of material during construction. In the optimization process, we have to ensure that only one material type is allocated at a specific supply location.

Mathematically, a set of binary decision variables $X_{i,l}$ is defined and controlled by constraint sets (25) and (26). In Eq. (25), for each material type $l = \{1, L\}$, where L is the total number of material types to be considered, there must be one assigned supply location within a site. Similarly, for each supply location $i = \{1, I\}$ where I is the total number of available supply points in a site that can store the construction material, at most one material type can be allocated as given in Eq. (26).

$$\sum_{i=1}^{I} X_{i,l} = 1, \qquad \forall l \in \{1, L\} \tag{25}$$

$$\sum_{l=1}^{L} X_{i,l} \leq 1, \qquad \forall i \in \{1, I\} \tag{26}$$

Objective function is defined as the total material transportation costs for optimization and these costs depend on the actual amount of material flows associating with different supply and demand locations. A set of auxiliary binary variables $\delta_{i,j,k,l}$ is thus defined to represent such existence of material flows which equals to "1" if material type l at supply point i is transported by a tower crane at location k to a demand point j and "0" otherwise. With the constraint set (13), the decision variables $\chi_{i,l}$, $\forall_{i,l}$ represent the linkage between material l and supply location i, $\Delta_j$ , $\forall j$ represent demand location j and $\zeta_k$ , $\forall k$ represent the selection of the tower crane $k$th location. Numerically, if all $\chi_{i,l} = 1$, $\Delta_j = 1$ and $\zeta_k = 1$, then the linkage of material flow is established giving $\delta_{i,j,k,l} = 1$. And $\delta_{i,j,k,l} = 1$ for all other cases so that no transportation cost will be counted in the objective function.

$$M\left(1 - \delta_{i,j,k,l}\right) \geq \left(3 - X_{i,l} - \Delta_j - \zeta_k\right) \geq \left(1 - \delta_{i,j,k,l}\right), \forall i \in \{1, I\}, \forall j \in \{1, J\}, \forall k \quad (27)$$
$$\in \{1, K\}, \forall l \in \{1, L\}$$

With the auxiliary variable $\delta_{i,j,k,l}$ expressing the existence of material flows, the total cost for material transportation from various supply points to demand points by a tower crane can be calculated using Eq. (28) that can be set as an objective function for optimization in the present formulation. The total cost $TC^h$ is simply defined as the sum of all transportation costs between supply and demand locations by a tower crane located at location k according to the set of material flow variables $\delta_{i,j,k,l}$ for a homogenous material supply system. In Eq. (28), $Q_{l,j}$ is the required quantity of material l at demand point j, C is the cost per unit time in operating a tower crane, and $T_{i,j}^k$ is the actual transport time between supply location i and demand location j by a tower crane at location k. The total cost can be evaluated and set as an objective function for optimization in the present formulation.

$$TC^h = \sum_{i=1}^{I} \sum_{j=1}^{J} \sum_{l=1}^{L} \delta_{i,j,k,l} T_{i,j}^k Q_{l,j} C, \forall k \in \{1, K\} \quad (28)$$

Optimization of the tower crane position in a homogeneous material supply system can be formulated to optimize the objective function in Eq. (28) subject to constraint sets in Eqs. (15)–(27).

### 3.2. Non-homogeneous material supply without area size constraint scenario

For larger construction sites in size without much physical size restrictions, the actual space allocated for material storage areas can be relatively increased. Different material types are able to be stored at one location (usually requiring larger space). For the case of this non-homogeneous material supply system, the tower crane operations and movements can be refined to save crane movement times and costs. To ensure that each potential material demand location is served, there must be at least one supply point for the required material types. Mathematically, this can be controlled by a set of binary variable $y_{i,j}$ in the constraint set (29). In this set, for each actual demand location $j \in \{1, J\}$ where J is the total number of demand location whenever $\Delta_j = 1$, at least one supply location must be selected and even one supply point can simultaneously serve for numerous demand locations.

$$\sum_{i=1}^{I} y_{i,j} \geq \Delta_j, \qquad \forall j \in \{1, J\} \quad (29)$$

An auxiliary binary-type variable $\mu_{i,j,k}$ is newly defined and governed by the constraint set (30) to establish the links among different supply and demand locations and the tower crane position. Numerically, whenever $y_{i,j} = 1$, meaning that a supply location i is linked with demand location j, and $\zeta_k = 1$ with tower crane location k is selected, then the linkage must be established that is $\mu_{i,j,k} = 1$ so that the transportation cost can be calculated.

$$M\left(1 - \mu_{i,j,k}\right) \geq \left(2 - y_{i,j} - \zeta_k\right) \geq \left(1 - \mu_{i,j,k}\right), \forall i \in \{1, I\}, \forall j \in \{1, J\}, \forall k \in \{1, K\} \qquad (30)$$

The total transportation cost for non-homogeneous supply point $TC^n$ is given by Eq. (31) according to the established material flow linkage $\mu_{i,j,k}$. In Eq. (31), $Q_{l,j}$ is the required quantity of material type l at a demand point j. C is the unit time cost of operating a tower crane, and $T_{i,j}^k$ is the actual transport time between supply location i and demand location j by a tower crane located at position k.

Optimizing of the tower crane position in a non-homogeneous material supply system can be formulated to optimize the objective function in Eq. (31) subjected to constraint sets in Eqs. (15)–(24), and (29), (30).

$$TC^n = \sum_{i=1}^{I} \sum_{j=1}^{J} \mu_{i,j,k} T_{i,j}^k Q_{l,j} C, \forall k \in \{1, K\} \qquad (31)$$

### 3.3 Non-homogeneous material supply location with physical size constraint

For the construction sites in the urban areas, work space is very limited and the material storage areas are comparatively small. In that sense, each material supply point can only supply construction materials for one demand point within a construction site. Similar to the previous non-homogeneous material supply strategy, different materials can still be stored at one material supply location. Mathematically, the binary variable $y_{i,j}$ , which is used for identifying linkages of material supply and demand locations as given in Eq. (29), is governed by two additional constraint sets in Eqs. (32) and (33). In Eq. (32), for each demand location $j \in \{1, J\}$ where j is the total number of demand locations to be considered, it must be assigned one supply location to store the materials. Similarly, for each supply location $i \in \{1, I\}$ with i being the total number of available supply location in a site that can store the construction material. Due to the storage area restriction, each supply location can only allocate materials for one demand location as given in Eq. (33).

$$\sum_{i=1}^{I} y_{i,j} = 1, \qquad \forall j \in \{1, J\} \qquad (32)$$

and

$$\sum_{j=1}^{J} y_{i,j} = 1, \qquad \forall i \in \{1, I\} \qquad (33)$$

To optimize the total material transportation cost $TC^n$ with these storage area constraints, the objective function in Eq. (31) can be applied and the problem is subjected to constraint sets as described in Eqs. (15)–(24), (29), (30), and (32), (33).

## 4. NUMERICAL EXAMPLES

Case studies that are used by Huang et al. [1] resolved by CBO, ECBO, and VPS to compare the applicability and performance of them by previous studies. Numerical examples consist of modeling a material supply and demand system considering 3 material types, 9 available material supply locations and 9 demand locations in a site which also provides 12 possible locations to set up and operate a tower crane. Hoisting velocity of the hook $V_h$=60 m/min, the radial velocity $V_a$=53.3 m/min, and the slewing velocity of the tower crane brachial $V_w$=7.57 rad/min. The operating cost of a tower crane per unit of time C is assumed to be $1.92 cost unit per minute and the quantities of material demand $Q_{l,j}$ are 10 units for material type l=1, 20 units for material type l=2, and 30 units for material type l=3 for all the demand points. The parameter β indicates the degree of coordination of hook movement in vertical and horizontal planes during practical operation is taken to be 0.25, and the α specifies the degree of coordination of hook movement in radial and tangential directions in the horizontal plane is assumed to be 1.0 [10,11]. For demonstration purpose, all $\gamma_k$=1.0 assuming no significant differences among the available locations for the tower crane operation. Table 1 lists all the three-dimensional (x, y, z) coordinates of all the potential demand points, all the potential locations for the material supply points, and tower crane potential locations.

Table 1: Coordinates of the potential locations

| # | | 1 | 2 | 3 | 4 | 5 | 6 | 7 | 8 | 9 | 10 | 11 | 12 |
|---|---|---|---|---|---|---|---|---|---|---|---|---|---|
| Demand point j | X | 34 | 34 | 51 | 60 | 76 | 76 | 60 | 51 | 43 | | | |
| | Y | 41 | 51 | 65 | 65 | 51 | 41 | 26 | 25 | 44 | | | |
| | Z | 15 | 15 | 15 | 15 | 15 | 15 | 15 | 15 | 15 | | | |
| Supply point i | X | 73 | 83 | 87 | 73 | 55 | 35 | 22 | 36 | 55 | | | |
| | Y | 26 | 31 | 45 | 67 | 73 | 67 | 46 | 27 | 15 | | | |
| | Z | 2 | 2 | 1.5 | 1.5 | 1.5 | 0 | 0 | 1 | 1 | | | |
| Tower crane position, k | X | 45 | 65 | 65 | 45 | 51 | 60 | 70 | 70 | 60 | 51 | 42 | 42 |
| | Y | 36 | 36 | 57 | 57 | 33 | 33 | 41 | 52 | 58 | 58 | 52 | 41 |
| | Z | 30 | 30 | 30 | 30 | 30 | 30 | 30 | 30 | 30 | 30 | 30 | 30 |

## 5. RESULTS AND DISCUSSION

In this study 30 independent experimental runs are performed for each scenario through 200 iterations. Employing three optimization methods, the problem is solved by MATLAB 2013.a [23]. Since the performance of the ECBO, and VPS are dependent on the control parameters, several tests have been conducted to select the appropriate parameters for finite-time performance of these algorithms.

### 5.1 Results and discussion for homogeneous material supply point scenario

The results of previous researches and relevant outputs including the total costs optimized locations of the homogeneous supply points *i* for the three material types l and a tower crane location $k$ are illustrated in Table 2. By comparison, it can be found that all three methods which are used in this paper can achieve to the results obtained by the MILP approach and almost 7% less than that those of GA. From the Table 2 and Fig. 8, it can be seen that the standard deviation and mean cost of the results of VPS is smaller than ECBO, and by ECBO is smaller than CBO. In addition, according to Table 2, supply points 2, 5, and 1 are selected for material types 1, 2, and 3, respectively to move by tower crane located at 8.

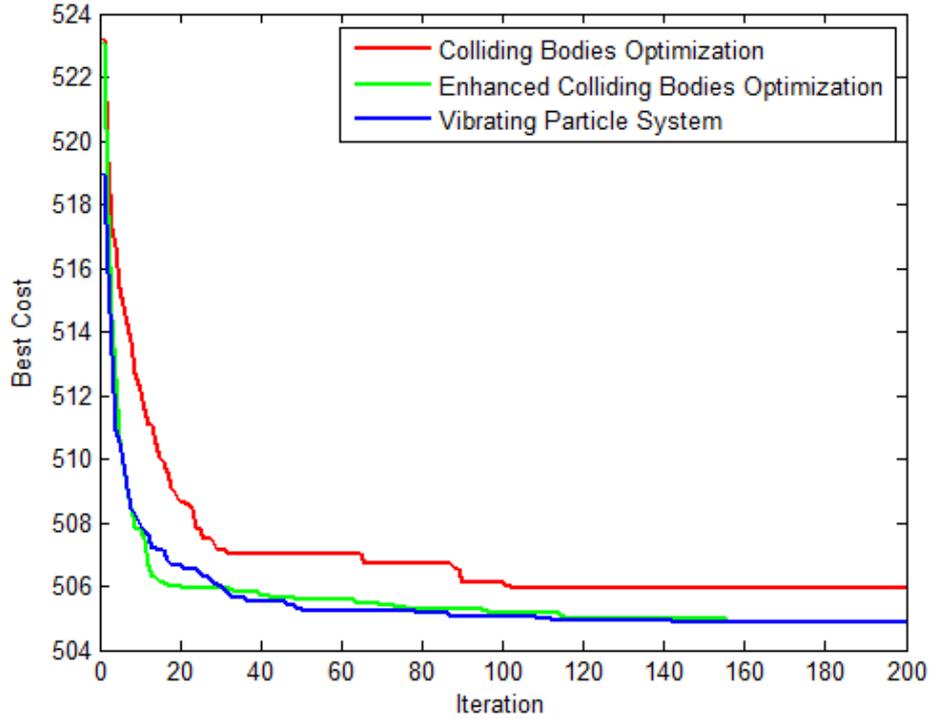

Figure 8. Mean cost of the homogeneous material supply point problem

Table 2: Comparison of the optimized design for homogeneous material supply point scenario

| method | Tower crane, k | Order of allocation of supply points to material type | | | Best cost | Mean cost | Standard deviation | Worst cost |
|---|---|---|---|---|---|---|---|---|
| | | 1 | 2 | 3 | | | | |
| GA (9) | 2 | 3 | 2 | 9 | 540.7587 | N/A | N/A | N/A |
| MILP (1) | 8 | 2 | 5 | 1 | 504.7631 | N/A | N/A | N/A |
| CBO | 8 | 2 | 5 | 1 | 504.7631 | 505.9426 | 1.3319 | 508.2809 |
| ECBO | 8 | 2 | 5 | 1 | 504.7631 | 504.8804 | 0.6423 | 508.2809 |
| VPS | 8 | 2 | 5 | 1 | 504.7631 | 504.8383 | 0.4121 | 507.0204 |

Note: N/A: Not available

### 5.2 Results and discussion for non-homogeneous material supply without area size constraint scenario

Table 3 shows the best costs and optimal design for non-homogeneous material supply without area size constraint scenario. By comparison, it can be found that all three methods used in this paper are upper than the results obtained by the MILP approach. And, as seen in Table 3 and Fig. 9, the mean total cost and standard deviation for ECBO is better than VPS, and for VPS is better than CBO. Thus, ECBO obtained a more stable evolution result than VPS and CBO.

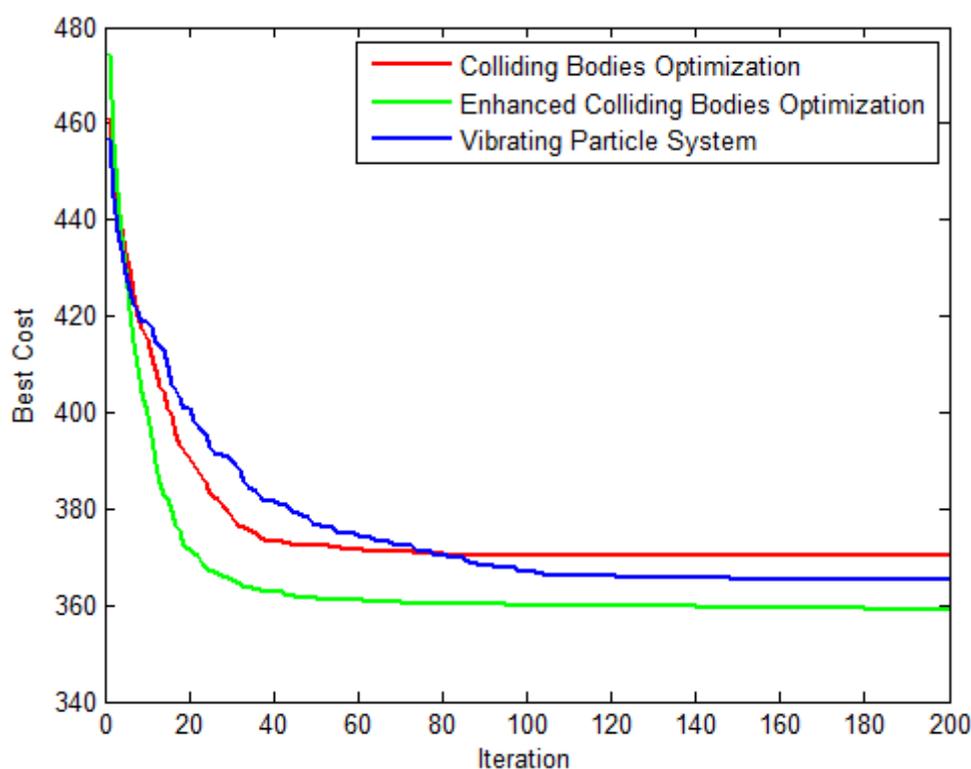

Figure 8. Mean cost of non-homogeneous material supply without area size constraint problem

Table 3: Comparison of the optimized design in non-homogeneous material supply without area size constraint scenario

| method | Tower crane, k | Order of supply point allocation | | | | | | | | | Best cost | Mean cost | Standard deviation | Worst cost |
|--------|------|---|---|---|---|---|---|---|---|---|------|------|------|------|
|        |      | 1 | 2 | 3 | 4 | 5 | 6 | 7 | 8 | 9 |      |      |      |      |
| MILP (1) | 9 | 3 | 7 | 4 | 4 | 3 | 2 | 1 | 1 | 3 | 343.3390 | N/A | N/A | N/A |
| CBO | 8 | 7 | 7 | 6 | 4 | 3 | 2 | 1 | 1 | 1 | 356.6403 | 370.4403 | 7.7215 | 391.2385 |
| ECBO | 8 | 7 | 7 | 6 | 4 | 3 | 2 | 1 | 1 | 1 | 356.6403 | 359.1270 | 4.4570 | 370.6018 |
| VPS | 8 | 7 | 7 | 6 | 4 | 3 | 2 | 1 | 1 | 1 | 356.6403 | 365.0899 | 8.2287 | 383.4908 |

Note: N/A: Not available;

*5.3 Results and discussion for non-homogeneous material supply with area size constraint scenario*

The best costs and optimal design for non-homogeneous material supply without area size constraint scenario are shown in Table 4. By comparison, it can be found that all of the utilized methods in this paper attained the results obtained by the MILP approach. And, as can be seen from Table 4 and Fig. 9, the mean total cost and standard deviation for the ECBO is better than others. Furthermore, ECBO obtained a more stable evolution result than VPS and CBO. Similar to MILP results, tower crane position 2 is selected and supply point locations are selected in the order of 7, 6, 5, 4, 3, 2, 1, 9, and 8 for demand points 1 thorough 9 in the best result of both methods.

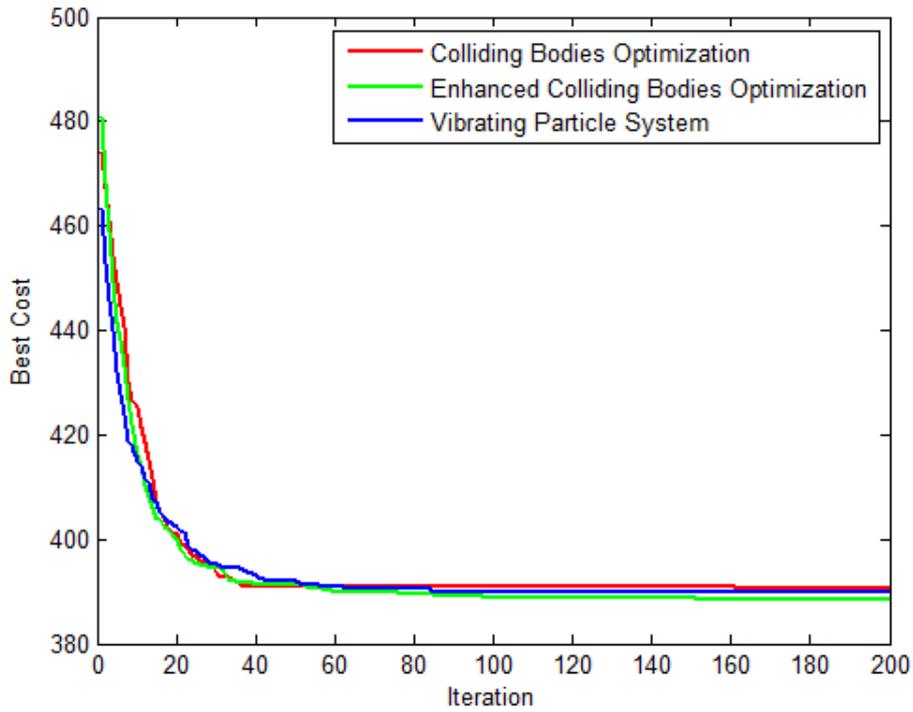

Figure 9. Mean cost of non-homogeneous material supply with area size constraint problem

Table 4: Comparison of the optimized results of non-homogeneous material supply for area size constraint

| method | Tower crane, k | Supply point allocation order to demand point | | | | | | | | | Best cost | Mean cost | Standard deviation | Worst cost |
|--------|------|---|---|---|---|---|---|---|---|---|-----------|-----------|-----------|-----------|
| | | 1 | 2 | 3 | 4 | 5 | 6 | 7 | 8 | 9 | | | | |
| MILP (1) | 2 | 7 | 6 | 5 | 4 | 3 | 2 | 1 | 9 | 8 | 388.2046 | N/A | N/A | N/A |
| CBO | 2 | 7 | 6 | 5 | 4 | 3 | 2 | 1 | 9 | 8 | 388.2046 | 390.649 | 2.7755 | 396.5171 |
| ECBO | 2 | 7 | 6 | 5 | 4 | 3 | 2 | 1 | 9 | 8 | 388.2046 | 388.3576 | 0.5940 | 391.3489 |
| VPS | 2 | 7 | 6 | 5 | 4 | 3 | 2 | 1 | 9 | 8 | 388.2046 | 389.9462 | 2.5784 | 398.9463 |

Note: N/A: Not available

## 6. CONCLUSIONS

In this paper, three newly developed meta-heuristic methods are employed for tower crane and material supply locations problem. Results show that except for homogeneous material supply point problem, ECBO presents more stable solution than VPS for all of considered scenarios. Also, both of ECBO, and VPS presents better solutions than CBO. However, the solution of this study for non-homogeneous material supply without area size constraint scenario could not reach to the solution obtained by Ref. [1].